\documentclass{elsarticle}
\biboptions{longnamesfirst,semicolon}
\usepackage{amsmath, amsthm, amssymb, latexsym}
\usepackage{graphicx}
\usepackage{amsfonts}
\usepackage{eucal}
\usepackage{amscd}
\usepackage{amssymb}
\usepackage{xypic}

\author[Il]{Paul Spiegelhalter \fnref{fn1}}
\ead{spiegel3@illinois.edu}

\author[Il]{Joseph Vandehey\corref{cor} \fnref{fn2}}
\ead{vandehe2@illinois.edu}

\cortext[cor]{Corresponding author}
\fntext[fn1]{Email address: spiegel3@illinois.edu}
\fntext[fn2]{Email address: vandehe2@illinois.edu}

\address[Il]{Department of Mathematics, 250 Altgeld Hall, 1409 W Green St., Urbana, IL, 61801}

\title{Squares in polynomial product sequences}
\date{\today}

\newtheorem{thm}{Theorem}[section]
\newtheorem{cor}[thm]{Corollary}
\newtheorem{lem}[thm]{Lemma}
\newdefinition{rem}{Remark}

\newtheorem{prop}[thm]{Proposition}
\newproof{pf}{Proof}
\newproof{pot}{Proof of Theorem \ref{thm2}}
\newcommand{\Mod}[1]{{\mbox{ (mod }#1 \mbox{)}}}

\begin{document}

\begin{abstract}
Let $F(n)$ be a polynomial of degree at least 2 with integer coefficients.  We consider the products $N_x=\prod_{1 \le n \le x} F(n)$ and show that $N_x$ should only rarely be a perfect power.  In particular, the number of $x \le X$ for which $N_x$ is a perfect power is $O(X^c)$ for some explicit $c<1$.  For certain $F(n)$ we also prove that for only finitely many $x$ will $N_x$ be squarefull and, in the case of monic irreducible quadratic $F(n)$, provide an explicit bound on the largest $x$ for which $N_x$ is squarefull.
\end{abstract}

\maketitle

\section{Introduction}

Several papers have recently been published concerning how often $$ N_x =\prod_{n\le x} F(n) $$ can be a perfect square, given an irreducible polynomial $F(n)$ with integer coefficients.  Cilleruelo proved in \cite{Cil} that if $F(n) = n^2+1$ then $N_x$ is a perfect square only when $x=3$.  Fang, using Cilleruelo's method, proved in \cite{Fang} that if $F(n) = 4n^2 +1$ or $F(n) = 2n(n-1) + 1$ then $N_x$ is never a perfect square, and G\"{u}rel and \"{O}zg\"{u}r Ki\c{s}isel proved in \cite{G} that if $F(n)=n^3+1$ then $N_x$ is never squarefull.  Conjectures regarding these products were initially put forth in \cite{Amd}, as they related to studying arithmetical properties of the arctangent function.

Cilleruelo, et al., in \cite{Cil2}, later showed that if $F(n)$ is an irreducible polynomial of degree at least $2$, then the number of times $N_x/d$ is a perfect square for $x$ in the interval $[M,M+N]$ is $$\ll N^{11/12} \left(\log{N}\right)^{1/3}$$ uniformly over all positive square-free integers $d$ and all positive integers $M$.

In this paper, we examine how often $ N_x =\prod_{n\le x} F(n) $ will be a perfect power or squarefull for more general $F(n)$.

If $F(n)$ is an irreducible monic quadratic, then we can provide an explicit bound on the largest $x$ for which $N_x$ can be squarefull.

We will also show that if $F(n)$ can be factored into linear and quadratic terms and given some conditions on the leading coefficient of the linear terms and discriminant of the quadratic terms, then $N_x$ will be squarefull for only finitely many $x$.  These conditions are general enough to cover some large collections of polynomials $F(n)$, such as all $F(n)$ that are the product of two or three distinct irreducible quadratics.  However, these proofs are not strong enough to provide explicit bounds.

More generally, we can show that if $F(n)$ is not of the form $s G(n)^p$ where $s$ is a rational number and $G(n)\in \mathbb{Z}[n]$, then $N_x$ is a perfect $p^{th}$ power for at most $O(X^{c_p})$ of the $x < X$ for some explicit $c_p < 1$.

In this paper, all polynomials denoted by lower-case letters are assumed to be irreducible over the rationals and have integer coefficients.  We denote the discriminant of a quadratic polynomial $f_i(n)$ by $D_i$.  Also, we assume $x \ge 1$ is integer-valued.

\section{The case $F(n) = n^2 +D$}
We wish to find an upper bound on those $x > 0$ for which
$$N_x := \prod_{n \leq x} \left(n^2 + D \right)$$
is squarefull.  Here, $D$ is a positive integer.  In particular, we will show that the bound $e^{C \cdot D}$ works, where $C$ is a constant that is effectively computable.
We start with the following proposition.
\begin{prop}\label{p:mainlemma}
If $a$ and $q$ are coprime natural numbers and $z$ a positive real number, then
$$S(z;q,a) := \left| \sum_{\substack{p \leq z \\ p \equiv a (q)}} \frac{\log p}{p} - \frac{1}{\phi(q)}\log z \right| = O(1),$$
where the constant implied by the Big-Oh expression can be effectively computed and is independent of $a$ and $q$.
\end{prop}

\begin{pf}
Here we use a method of proof similar to that employed by Pomerance in \cite{PomAmicable}.

Suppose that $z \geq e^{q^{2/3}}$.  Define
$$\theta(z;q,a) := \sum_{\substack{p \leq z \\ p \equiv a (q)}} \log p.$$
Now,
\begin{eqnarray*}
\sum_{\substack{p \leq z \\ p \equiv a (q)}} \frac{\log p}{p} &=& \frac{1}{z} \theta(z;q,a) - \frac{1}{2} \theta(z;q,a) + \int_2^z \frac{\theta(t;q,a)}{t^2} dt \\
&\leq& \frac{1}{z} \theta(z;q,a) + \left( \int_2^{q} + \int_{q}^{e^{q^{1/2}}} + \int_{e^{q^{1/2}}}^{e^{q^{2/3}}} + \int_{e^{q^{2/3}}}^z \right) \frac{\theta(t;q,a)}{t^2} dt. \\
\end{eqnarray*}
To bound the first term, we use the bound $\theta(z) \leq 2z\log 2$.  This follows from the inequality $\prod_{p \leq n} p \leq 4^n$ (see for example \cite{Rose}).  Hence
$$\frac{1}{z} \theta(z;q,a) \leq \frac{1}{z}\theta(z) \leq 2\log 2.$$
For the first of the four integrals, we note that each of
$$\sum_{\substack{p \leq q \\ p \equiv a (q)}} \frac{\log p}{p} \text{  and  } \frac{1}{q} \sum_{\substack{p \leq q \\ p \equiv a (q)}} \log p $$
is bounded by $1$, so by partial summation
$$\int_2^{q} \frac{\theta(t;q,a)}{t^2} dt \leq 3.$$
Since $\theta(z;a,q) \leq (1+z/q) \log{z} \leq \frac{2z \log z}{q}$ when $z \geq q$,
$$ \int_{q}^{e^{q^{1/2}}} \frac{\theta(t;q,a)}{t^2} dt \leq \int_{q}^{e^{q^{1/2}}} \frac{2\log t}{q t} dt \leq 1.$$
Now, the Brun-Titchmarsh theorem in the form of Montgomery and Vaughan (see \cite{MontgomeryVaughan}) gives us that
$$\pi(z;q,a) \leq \frac{2z}{\phi(q)\log (z/q)}$$
for $z>q$.  So for $z> e^{q^{1/2}}$,
$$\theta(z;a,q) \leq \frac{2z}{\phi(q)}\left( \frac{1}{1 - \frac{\log q}{\log z}} \right) \leq \frac{8z}{\phi(q)},$$
using elementary calculus.  Hence
$$\int_{e^{q^{1/2}}}^{e^{q^{2/3}}} \frac{\theta(t;q,a)}{t^2}dt \leq \frac{8}{\phi(q)} \int_{e^{q^{1/2}}}^{e^{q^{2/3}}} \frac{dt}{t} \leq \frac{8q^{2/3}}{\phi(q)}.$$
If $z \geq e^{q^{2/3}}$ then by the prime number theorem for arithmetic progressions (see \cite{Davenport} p. 123),
$$\left| \theta(z;q,a) - \frac{z}{\phi(q)} \right| \leq Aze^{-c(\log z )^{1/8}} $$
where $A$ and $c$ are positive absolute constants.  Then since
$$\int_{e^{q^{2/3}}}^z \frac{1}{\phi(q) t} dt = \frac{1}{\phi(q)} \log z - \frac{q^{2/3}}{\phi(q)},$$
and since $Ae^{-c(\log z )^{1/8}} \leq \frac{A'}{(\log z)^2}$ for some $A'$ depending on $A$ and $c$, we have
\begin{eqnarray*}
\left| \int_{e^{q^{2/3}}}^z \frac{\theta(t;q,a)}{t^2} dt - \frac{1}{\phi(q)}\log z \right| &\leq& \frac{q^{2/3}}{\phi(q)} + \int_{e^{q^{2/3}}}^z \frac{Ae^{-c(\log t)^{1/2}}}{t} dt \\
&\leq& \frac{q^{2/3}}{\phi(q)} + \int_{e^{q^{2/3}}}^z \frac{A'}{t \log^2 t} dt \\
&\leq& \frac{q^{2/3}}{\phi(q)} + \frac{A'}{q^{2/3}}
\end{eqnarray*}
Thus
\begin{eqnarray*}
S(x;q;a) &=& \left| \sum_{\substack{p \leq z \\ p \equiv a (q)}} \frac{\log p}{p} - \frac{1}{\phi(q)}\log z \right| \\
&\leq& 4 + 2 \log 2 + \frac{9q^{2/3}}{\phi(q)} + \frac{A'}{\phi(q)}.
\end{eqnarray*}
Since $q^{2/3}/\phi(q)$ can be effectively bounded, this completes the proof of the proposition for $z > e^{q^{2/3}}$.  For smaller values of $z$, one can simply truncate the expansion of
$$\sum_{\substack{p \leq z \\ p \equiv a (q)}} \frac{\log p}{p}$$
as a sum of integrals at the appropriate place to obtain a similar bound. $\Box$ \end{pf}

\begin{rem}
Using the Brun-Titchmarsh estimate we can show that
$$\frac{1}{z} \theta(z;q,a)  = \frac{\log{p(q;a)}}{p(q;a)} + O\left(\frac{\log q}{q}\right)$$
where $p(q;a)$ denotes the first prime $p \equiv a \Mod q$ and that, after a suitable adjustment to the bounds of integration, the remaining terms are $O(q^{-1/3})$.  One obtains the result
$$S(z;q,a) = \frac{\log{p(q;a)}}{p(q;a)} + O(q^{-1/3})$$
for $x>e^{q^{2/3}}$, where the implied constant is independent of $a$ and $q$.
\end{rem}

We have that if $E$ is a set of residue classes mod $q$, then
$$\left| \sum_{\substack{p \leq x \\ p \in E}} \left( \frac{\log p}{p-1} - \frac{\log p}{p} \right) \right| \leq 1. $$
Thus as a corollary to Proposition \ref{p:mainlemma} we have that for some constant $C_0$
\begin{equation}\label{eq:firstcorollary}
\sum_{\substack{p \leq x \\ p \in E}} \frac{\log p}{p-1} \leq \frac{|E|}{\phi(q)} \log x +  |E|C_0
\end{equation}

\begin{prop}\label{p:littleLemma}
The number $N_x$ satisfies
$$\log N_x \geq 2x \log x - 2x.$$
\end{prop}
\begin{pf}
Note that $e^x \geq \frac{x^x}{x!}$ by Taylor series, so $x! \geq \left( \frac{x}{e} \right)^x$, hence
$$\sum_{n \leq x} \log{n} \geq x\log{x} - x.$$
Then since $\log(x^2 + D) \geq 2 \log x$, we have
$$\log N_x \geq 2 \sum_{n \leq x} \log n  \geq  2x \log x - 2x.$$ $\Box$ \end{pf}

\begin{prop}
There is a prime factor $p_x$ of $N_x$ satisfying $p_x> \frac{1}{72}x\log x$ for all $x$ larger than $C_1e^{C_2 D} = \exp \{(\frac{8}{5}((4C_0+8)D+2)\}$, where $C_0$ is the constant defined in ~(\ref{eq:firstcorollary}).
\end{prop}
\begin{pf}
Let $k = \frac{1}{72}$.  For a given $x$, let $\alpha_p$ be defined for each prime $p$ so that $N_x = \prod_p p^{\alpha_p}$.  Now, $p|N_x$ only when $p|D$, or when $p\nmid D$ and $-D$ is a quadratic residue mod $p$.  The latter occurs only for a particular set $S$ of residue classes mod $4D$ with $2\left| S \right| = \phi(4D)$.  Hence
$$N_x = \prod_{p|D} p^{\alpha_p} \prod_{p \in S} p^{\alpha_p}$$
where by a slight abuse of notation we take $p \in S$ to mean $p \Mod {4D} \in S$.  Now, if $p \nmid D$, then each interval of length $p^j$ contains at most $2$ solutions of $n^2+D \equiv 0 \Mod {p^j}$.  So
\begin{eqnarray}
\alpha_p &\leq& \sum_{j \leq \frac{\log (x^2 + D)}{\log p}} 2 \lceil x/p^j \rceil \label{eq:alphaNotD}\\
&\leq& 2x \sum_{j \leq \frac{\log (x^2+D)}{\log p}} \frac{1}{p^j} + 2\frac{\log (x^2 + D)}{\log p } \nonumber \\
&\leq& \frac{2x}{p-1} + 2\frac{\log (x^2 + D)}{\log p } \nonumber
\end{eqnarray}
On the other hand, if $p|D$, we write $D = p^{e_0}D'$, with $p \nmid D'$.  Then as a result of Huxley (see \cite{Huxley}) we have that $n^2 + D \equiv 0 \Mod {p^j}$ has at most $2p^{e_0}$ solutions.  By an argument similar to that in~(\ref{eq:alphaNotD}) we have
\begin{equation}\label{eq:alphaD}
\alpha_p \leq \frac{2p^{e_0}x}{p-1} + 2p^{e_0}\frac{\log (x^2 + D)}{\log p }.
\end{equation}
If the claim in the proposition does not hold, then there is an $x > C_1 e^{C_2 D}$ such that
$$N_x = \prod_{\substack{p \leq kx \log x \\ p|D}} p^{\alpha_p} \prod_{\substack{p\leq kx \log x \\ p \in S}} p^{\alpha_p}.$$
To estimate $\log N_x$ with $N_x$ in this form, we use Chebyshev's inequality $\pi(x) < 2 \frac{x}{\log x}$ as given in \cite{CrandalPomerance}.  Also, note that for $x$ in the prescribed range
\begin{equation}\label{eq:log}
\log(x^2 + D) \leq 3 \log x.
\end{equation}
Since $k > x^{-1/3}$ we have $\log(kx\log x) \geq \frac{2}{3}\log x$, so
\begin{equation}\label{eq:piPoly}
\pi(kx \log x) \log( x^2 + D)\leq 6 \frac{k x \log^2 x}{\log(kx \log x )} \leq \frac{1}{8}x \log x.
\end{equation}
Now, certainly we have that $kx\log x > D$ for $x > C_1 e^{C_2 D}$, so by ~(\ref{eq:alphaD}) and ~(\ref{eq:log})
\begin{eqnarray*}
\sum_{\substack{p \leq kx \log x \\ p|D}} \alpha_p \log p &\leq& 2x \sum_{p|D}p^{e_0} \frac{\log p}{p-1} + 2  \sum_{p|D} p^{e_0} \frac{\log (x^2 + D)}{\log p }\\
&\leq& 2x \sum_{p|D}p^{e_0} + 2 \log(x^2 + D) \sum_{p|D}p^{e_0} \\
&\leq& 2x \prod_{p|D}p^{e_0} + 2 \left( 3 \log{x} \right) \prod_{p|D}p^{e_0} \\
&\leq& 8Dx.
\end{eqnarray*}
Now, if $x$ is in the prescribed range then $\log x \leq x^{1/8}$, so by ~(\ref{eq:alphaNotD}), ~(\ref{eq:piPoly}) and ~(\ref{eq:firstcorollary}) we have
\begin{eqnarray*}
\sum_{\substack{p \leq kx \log x \\ p \in S}} \alpha_p \log p &\leq& 2x \sum_{\substack{p \leq kx \log x \\ p \in S}} \frac{\log p}{p-1} + 2 \log(x^2+D) \sum_{\substack{p \leq kx \log x \\ p \in S}} 1\\
&\leq& 2x \left( \frac{|S|}{\phi(4D)}\log(kx \log x) +|S|C_0 \right) + 2 \pi(kx \log x) \log(x^2+D)\\
&\leq& x \log(k x^{9/8}) + \phi(4D)C_0x + \frac{1}{4} x \log x \\
&\leq& \frac{11}{8} x \log x + 4 C_0 D x.
\end{eqnarray*}
So
\begin{eqnarray*}
\log N_x &=& \sum_{\substack{p \leq kx \log x \\ p|D}} \alpha_p \log p + \sum_{\substack{p \leq kx \log x \\ p \in S}} \alpha_p \log p\\
&\leq& \frac{11}{8} x \log x + 4 C_0 D x + 8 D x
\end{eqnarray*}
Thus by Proposition \ref{p:mainlemma} we have
$$\frac{5}{8} x \log x \leq (4 C_0 D +8D + 2) x $$
hence $x \leq e^{8((4C_0 + 8)D + 2)/5}$.  This is a contradiction to our earlier assumtion that $x > C_1 e^{C_2 D}$. $\Box$ \end{pf}
\begin{thm}\label{th:theorem}
For any $x$ larger than $C_1 e^{C_2 D}$ the number $N_x$ is not squarefull.
\end{thm}
\begin{pf}
If $N_x$ is squarefull and $p|N_x$, then either $p^2|n^2+D$ for some positive integer $n \leq x$ or $p|n^2+D$ and $p|m^2+D$ for some distinct positive integers $n,m \leq x$.  In the first case, we have
$$p\leq \sqrt{x^2+D} \leq x + D \leq 2x$$
since $x>D$.  In the second case, we have that $p$ divides $n^2 - m^2 = (n-m)(n+m)$, so $p \leq 2x$.  If $x$ is in the range given in the theorem, then
$$2x < k x \log x < p_x$$
for some $p_x$ dividing $N_x$, a contradiction. $\Box$ \end{pf}

\begin{rem}
If we take $F(n)$ to be any irreducible monic quadratic, we can apply the above technique to $|N_x|$ to obtain similar results.  Write $F(n) = (n - \alpha)^2 +D$; there is a constant $C_f < 0$ such that
$$\log \left|1- \frac{2 \alpha}{n} + \frac{\alpha^2 + D}{n^2} \right| > C_f.$$
Modifying Proposition ~\ref{p:littleLemma} we get
\begin{eqnarray*}
\log |N_x| &\geq& 2 \sum_{n \leq x} \log{n} + C_f x \\
&\geq& 2x \log{x} - 2x + C_f x.
\end{eqnarray*}
The rest of the proof of Theorem ~\ref{th:theorem} holds with only slight modification.  One obtains the result that $N_x$ is not squarefull for any $x$ larger than $\exp\{\frac{8}{5}(2 - C_f +4C_0 D + 8D) \}$.
\end{rem}

\section{Products of quadratics}

The main result of this section relies on the following theorem, proved in two separate cases by Duke, Friedlander, and Iwaniec in \cite{DFI} and T{\'o}th in \cite{Toth}.

\begin{thm}\label{DFIT}

If  $f(n)$ be an irreducible quadratic polynomial with integer coefficients, and $0 \le \alpha < \beta \le 1$, then $$K_x= \# \left\{ (p,v) | 0 \le v < p \le x, f(v) \equiv 0 \pmod{p},\alpha \le \frac{v}{p} < \beta \right\} \sim (\beta - \alpha)\pi (x)$$ where $v \in \mathbb{Z}$, $p$ prime, and the asymptotic relation holds as $x \to \infty$.

\end{thm}

We begin by presenting two lemmas derived from this result, which we will often refer to as the DFIT result, after its various authors.

\begin{lem}\label{DFITlem}
Let $f(n)=n^2+bn+c$ be a monic quadratic polynomial, and let $\epsilon > 0$.  Then there exists $\delta=\delta(\epsilon)$ and $x_0$ such that for all $x > x_0$, at least $(\frac{1}{2}-\epsilon)(\pi(2x)-\pi((2-\delta)x))$ of the primes between $(2-\delta)x$ and $2x$ divide $N_x=\prod_{n\le x} f(n)$ exactly one time.

\end{lem}

\begin{pf}

In particular we will choose $\delta$ such that $1/2+\epsilon /2 > 1/(2-\delta)$.

We first note that for all sufficiently large $x$, $(2-\delta)^2 x^2 > f(x)$, so any prime $p \in \left[(2-\delta)x, 2x\right]$ can only divide any given $f(n)$ at most one time.  Thus the number of times $p$ divides $N_x$ equals the number of $n \le x$ for which $f(n) \equiv 0 \pmod{p}$.

We rewrite this last condition as $$\left( n+\frac{b}{2} \right)^2 \equiv \frac{b^2 - 4 c}{4} \pmod{p}$$ and then write $b/2$ in reduced terms as $B/A$ and first handle the case where $A=1$.

Then if we use the interval $$\left( \frac{1}{2-\delta}+\frac{\epsilon}{4} , 1-\frac{\epsilon}{4} \right)$$ and the polynomial $\overline{f}(n)=(2n)^2-(b^2-4c)$ in the DFIT result, we see that the number of pairs $(v,p)$, for which $0 \le v < p$, $(2-\delta)x < p < 2x$, $p | \overline{f}(v)$ and $$\frac{v}{p} \in \left( \frac{1}{2-\delta}+\frac{\epsilon}{4}, 1-\frac{\epsilon}{4} \right),$$ tends asymptotically to $\delta (1-\frac{\epsilon}{2} - \frac{1}{2-\delta}) \frac{x}{\log{x}}$, as $K_{2x} \sim 2 (1-\frac{\epsilon}{2} - \frac{1}{2-\delta}) \frac{x}{\log{x}}$ and $K_{(2-\delta)x} \sim (2-\delta) (1-\frac{\epsilon}{2} - \frac{1}{2-\delta}) \frac{x}{\log{x}}$.

If $\overline{f}(v) \equiv 0 (\textrm{mod} \quad p)$, then $v^2 \equiv (b^2 - 4 c)/4$.  If we set $n=v-B$, this gives us a solution to $f(n)\equiv 0 (\textrm{mod} \quad p)$.   We can pick $x$ to be large enough so that $B/p < \epsilon / 20$.  So, $$\frac{n}{p} \in  \left( \frac{1}{2-\delta}+\frac{\epsilon}{5}, 1-\frac{\epsilon}{5} \right)$$ with $0 < n < p$.  This implies that $$n > p \left(\frac{1}{2-\delta}+\frac{\epsilon}{5}\right)> (2-\delta)x \left(\frac{1}{2-\delta}+\frac{\epsilon}{5}\right) > x + \frac{\epsilon}{5}.$$  So this particular $p$ can only divide $N_x$ at most once.

Moreover, each pair $(v,p)$ corresponds in a one to one ratio with pairs $(p-v,p)$, with $0 \le p-v < p$, $(2-\delta)x < p < 2x$, $p | \overline{f}(p-v)$ and $$\frac{v}{p} \in \left( \frac{1}{2-\delta}+\frac{\epsilon}{4}, 1-\frac{\epsilon}{4} \right),$$ which is the same thing as $$\frac{p-v}{p} \in \left( \frac{\epsilon}{4}, 1-\frac{1}{2-\delta}-\frac{\epsilon}{4} \right).$$  Again setting $n=p-v-B$ and extending the bounds to allow $$\frac{n}{p} \in \left( \frac{\epsilon}{5}, 1-\frac{1}{2-\delta}-\frac{\epsilon}{5} \right),$$ we can see that $$n < p (1-\frac{1}{2-\delta}-\frac{\epsilon}{5}) < 2x (1-\frac{1}{2-\delta}-\frac{\epsilon}{5} ) < x,$$ so that this $p$ must divide $N_x$ at least once, and hence, by the last paragraph, exactly once.

As there are asymptotically $\delta \frac{x}{\log{x}}$ primes in the interval $(2-\delta)x $ to $ 2x$, and our choice of $\delta$ implies $$1-\frac{\epsilon}{2} - \frac{1}{2-\delta} > \frac{1}{2} - \epsilon,$$ we have proved the lemma in this case.

For the case $A=2$, we need to consider how $B/A$ acts modulo $p$.  For all odd primes, $1/2 \equiv (p+1)/2$.  Since $p$ is odd, $(p+1)/2$ is an integer so this represents a solution to $1/2 \pmod{p}$.  Therefore $B/2 \equiv B (p+1)/2$.  If we call this latter integer $k$, $0 \le k < p$, then note that $k/p$ tends towards $1/2$ as $p$ grows since $B$ is a fixed odd number.

From here, the proof of the second case proceeds identically to that of the first case, except that we use the interval $$\left( \frac{1}{2-\delta} -\frac{1}{2}+\frac{\epsilon}{4}, \frac{1}{2}-\frac{\epsilon}{4}  \right)$$ in the DFIT result and set $n=v-k$ or $n=p-v-k$ as appropriate. $\Box$ \end{pf}

\begin{rem}
Clearly the previous proof also works if $f(n)=an^2+bn+c$, where $a|b$ or $2a|b$.  In general though, the $b/2a$ term is only well-behaved over primes of a specific congruence class, and the DFIT result does not address the equidistribution of $v/p$ for primes $p$ of a specific congruence class, so we do not yet know how to extend the above lemma.
\end{rem}

\begin{lem}\label{oneprime}

Let $f(n)$ be an irreducible monic quadratic polynomial with integer coefficients, and $2 < a < b$.  Then for all sufficiently large $x$, there exists a prime $p$, $ax < p < bx$, such that $p | \prod_{n \le x} f(n)$.

\end{lem}

\begin{pf}

Consider pairs $(p,v)$ for which p divides $f(v)$, with $ax \le p \le bx$, and $0 \le v/p \le 1/b$. By DFIT, the number of such pairs is asymptotically $(1-a/b) x/\log{x}$. In particular, there is always such a pair once $x$ is sufficiently large. But for this pair, we have $v \le p/b \le x$, so that $p$ divides $f(v)$, which itself divides $\prod_{n \leq x} f(n).$ $\Box$ \end{pf}

Cilleruelo, in his proof, used the fact that if $N_x$ is a perfect square, then all primes dividing it must be less than $2x$, so the previous lemma provides an alternative proof that $\prod_{n \le x} (n^2+1)$ is not infinitely often a square.  We can generalize this idea a little further with the help of the following lemmas.

\begin{lem}

If $f_1(n) = a_1 n^2 + b_1 n + c_1$ and $f_2(n) = a_2 n^2 + b_2 n + c_2$ are two distinct quadratic polynomials such that $D_1 D_2$ is a square, then the largest prime $p$ that can divide both $N_x = \prod_{n \le x} f_1(n)$ and $M_x = \prod_{n \le x} f_2(n)$ is bounded by $cx$ for some positive constant $c$ and sufficiently large $x$

\end{lem}

\begin{pf}

We can rewrite $f_1(n) = a_1(n+(b_1/2a_1))^2 - (b_1^2 - 4a_1 c_1)/4 a_1$ and $f_2(n)= a_2(n+(b_2/2a_2))^2 - (b_2^2 - 4a_2 c_2)/4 a_2$.  Thus writing $d=\sqrt{D_1D_2}$, we have that \begin{align*}
&4a_1 D_2 f_1(n) - 4a_2 D_1 f_2(m) \\
&= 4a_1^2 D_2 \left(n+\frac{b_1}{2 a_1}\right)^2 -4a_2^2 D_1\left(m+\frac{b_2}{2 a_2}\right)^2\\
&= \frac{1}{D_1} \left( D_1D_2 \left(2a_1  \left(n+\frac{b_1}{2 a_1}\right)\right)^2 -\left(2a_2 D_1 \left(m+\frac{b_2}{2 a_2}\right)\right)^2 \right)\\
&=\frac{1}{D_1} \left( d  \left(2a_1n+b_1\right) - D_1\left( 2a_2m+b_2\right) \right) \left( d  \left(2a_1n+b_1\right) + D_1\left( 2a_2m+b_2\right) \right).
\end{align*} So if $p|f_1(n)$ and $p|f_2(m)$, then $p$ must divide the second or third factor of the above equation.  Since $n,m \le x$ by assumption, this implies that $p$ must be less than $\left|d  (2a_1x+b_1)\right| + \left|D_1( 2a_2x+b_2)\right| < (2|a_1 d| + 2|a_2 D_1| +1) x$ for sufficiently large $x$. $\Box$ \end{pf}

\begin{lem}

If $f(n) = n^2 + b n + c$ is a quadratic polynomial and for a prime $p$, $p^2 | \prod_{n \le x} f(n)$, then $p < (2+|b|+|c|)x$.

\end{lem}

\begin{pf}

If $p^2 | f(n)$ for some $n \le x$, then $p^2 \le n^2 + bn + c \le x^2 + |b|x+|c| < (1+|b|+|c|) x^2$, which implies that $p \le x\sqrt{1+|b|+|c|}$.

If $p | f(n)$ and $p | f(m)$ for some $n,m \le x$, then $p | n^2 + bn + c - m^2 - bm - c = (n-m)(n+m) + b(n-m) = (n-m)(n+m+b)$.  Since $p$ is prime, this implies $p | (n-m)$ or $p | (n+m+b)$.  Either way this implies that $p < (2+|b|)x$.  So the lemma holds. $\Box$ \end{pf}

Thus we have the following result using our variant method of Cilleruelo.

\begin{thm}\label{ThmWeak}

Let $f_i(n)$, $1 \le i \le I$, be some sequence of monic irreducible polynomials.  If $D_1 D_i$ is a perfect square for all $1 \le i \le I$, then $$N_x = \prod_{n \le x} \prod_{i=1}^I f_i(n)$$ cannot be squarefull for infinitely many $x$.

\end{thm}

\begin{pf}

First, suppose $N_x$ is squarefull, so for all primes $p$ such that $p$ divides $N_x$, $p^2|N_x$.

If $I=1$ then by the previous lemma, there exists some constant $c$ independent of our choice of $x$, for which $p < cx$ for all primes dividing $N_x$.

If $I > 1$, then for each prime $p | N_x$, either for some $i$, $p^2 | \prod_{n\le x} f_i(n)$, or else for some $i$ and $i'$, $p | \prod_{n\le x} f_i(n)$ and $p | \prod_{n \le x} f_{i'}(n)$.  Regardless of which case we fall into, the previous two lemmas tell us that there exists some constant $c$, dependent only on the $f_i$'s for which $p < cx$ for all sufficiently large $x$.

But again, Lemma \ref{oneprime} shows that $\prod_{n \le x} f_1(n)$ will eventually be divisible by at least one prime in the range $cx$ to $(c+1)x$.  Thus our assumption that $N_x$ could be squarefull for any of these large $x$ must be false. $\Box$ \end{pf}

We can replace the condition that requires $D_1 D_i$ to be a perfect square through the use of the following lemma.

\begin{lem}\label{congruence}

If $f_i(n)$, $1 \le i \le I$, is some sequence of distinct irreducible polynomials, with $$J_f :=  1+ \sum_{\substack{ \varnothing \neq J \subset \{1,2,3,\ldots,I\} \\ \prod_{j \in J} D_j \textrm{square}}}  (-1)^{|J\setminus \{1\}|} > 0,$$ then there exists some residue class $k$ modulo $\prod_{i=1}^I D_i$, such that all sufficiently large primes congruent to $k \pmod{\prod_{i=1}^I D_i}$ cannot divide any term of the form $f_i(n)$ for $1 < i \le I$, but will divide some term of the form $f_1(n)$.

\end{lem}

\begin{pf}

Once again, a given prime $p$ will divide $f_i(n) = a_i n^2 + b_i n + c_i$ if $$(n+\frac{b_i}{2a_i})^2 + \frac{4a_i c_i - b_i^2}{4a_i^2} \equiv 0 \pmod{p},$$ which makes sense provided $p$ is larger than $a_i$.

Thus $p$ will divide $f_i(n)$ for some $n$ if and only if $D_i$ is a quadratic residue modulo $p$.  To estimate the number of primes up to $z$, which can divide $f_1(n)$ for some $n$ but can never divide $f_i(m)$ for $i \neq 1$, we use the formula $$\sum_{D < p \le z} \left( \frac{1+\left( \frac{D_1}{p} \right)}{2} \prod_{2 \le i \le I} (-1) \frac{-1+\left( \frac{D_i}{p} \right)}{2} \right).$$

Here, $D$ is some constant larger than all the $D_i$.  But this sum equals \begin{align*}
&(-1)^{|I|-1} \frac{1}{2^{I}} \sum_{D < p \le z} \sum_{J \subset \{1,2,3,\ldots,I\}} (-1)^{|I|-|J\backslash \{1\}|-1} \left( \frac{\prod_{j \in J} D_j}{p} \right)\\
&= \frac{1}{2^{I}} \sum_{J \subset \{1,2,3,\ldots,I\}} \sum_{D < p \le z}  (-1)^{|J\backslash \{1\}|} \left( \frac{\prod_{j \in J} D_j}{p} \right) \end{align*}

If $\prod_{j \in J} D_j$ is a square, then $$\sum_{D < p \le z} \left( \frac{\prod_{j \in J} D_j}{p} \right) $$ will be asymptotic to $\pi (z).$  Otherwise, the sum will be $o(\pi (z))$ (in fact, it is $O(1)$).  Thus the sum above equals $$\frac{\pi (z)}{2^{I}} \left( 1 + \sum_{\substack{\varnothing \neq J \subset \{1,2,3,\ldots,I\} \\ \prod_{j \in J} D_j \textrm{square}}}  (-1)^{|J\backslash \{1\}|} + o(1)\right) $$ which will represent a non-trivial proportion of the primes provided $ J_f > 0$. $\Box$ \end{pf}

We can now combine this with Lemma \ref{DFITlem}, assuming $f_1$ is monic.  If we pick $\epsilon < 1/\phi (D)$, then Lemma \ref{DFITlem} says that for all sufficiently large $x$ there exists a prime congruent to $k \pmod{\prod_{i=1}^I D_i}$ that must divide $\prod_{n\le x} f_1(n)$ exactly once. And since it cannot divide $f_i(n)$ for $1 < i \le I$, we have proved the following theorem.

\begin{thm}\label{NotSquarefull}

Suppose that we have a set of $I$ distinct irreducible quadratic polynomials $f_i(n)=a_i n^2+b_i n+c_i$ with $f_1$ monic.  Furthermore, suppose that $J_f > 0$

Then for sufficiently large $x$ the number $$ N_x = \prod_{n \le x} f_1(n) $$ is not squarefull.  Moreover, $N_x$ cannot be made a squarefull by multiplying $N_x$ with terms of the form $f_i(n)$ with $i \neq 1$, $n \in \mathbb{N}_{>0}$.

\end{thm}

\begin{cor}\label{ThmStrong}

Suppose that we have a set of $I$ distinct irreducible quadratic polynomials $f_i$ with $f_1$ monic.  Furthermore, suppose that $J_f > 0.$

Then the number

$$ N_x = \prod_{n \le x} \prod_{i=1}^I f_i(n) $$
cannot be infinitely often a squarefull.

\end{cor}

While the conditions of the previous theorems have been somewhat complex, we can combine them to prove the following - much simpler - theorem.

\begin{cor}

Suppose we have $k$ distinct monic quadratic polynomials $f_i$.  Then

$$N_x = \prod_{n \le x} \prod_{1 \le i \le k} f_i(n)$$

cannot be infinitely often squarefull if $k=2$ or $3$.

\end{cor}

\begin{pf}

If any $D_i$ is a perfect square, then $f_i$ is reducible, so none of the $D_i$ can be a perfect square.

In the case $k=2$, we therefore have only two cases to consider: either $D_1D_2$ is a perfect square or it is not.

If $D_1D_2$ is a perfect square, then we apply Theorem \ref{ThmWeak}.

If $D_1D_2$ is not a perfect square, then we apply Theorem \ref{ThmStrong} as $J_f = 1$ in this case.

In the case $k=3$ we again have multiple sub-cases.

First, if no product of $D_1,D_2,D_3$ is ever a square, we may again apply Theorem \ref{ThmStrong} as $J_f = 1$ in this case.

Suppose that exactly one product of two of the discriminants is a square, and that the product of all three is not.  By reindexing we can let $D_2D_3$ be the square.  Then we again apply Theorem \ref{ThmStrong} as $J_f = 2$ in this case.

Note that it is impossible to have just two products of two discriminants being square, as if $D_1 D_2$ and $D_1 D_3$ are square, then so is $(D_1 D_2)(D_1 D_3)/D_1^2 = D_2 D_3$.

So suppose that all three products of two of the discriminants is a square, and that the product of all three is not.  Here we can apply Theorem \ref{ThmWeak}.

Suppose that the only square can be formed by multiplying all three discriminants together, i.e. $D_1 D_2 D_3$ is a square, then we apply Theorem \ref{ThmStrong} as $J_f = 2$ in this case.

Finally assume that some product of two discriminants and the product of all three discriminants are squares, say $D_1 D_2$ and $D_1D_2D_3$ are both squares.  Then $(D_1D_2D_3)/(D_1 D_2)=D_3$ must also be a square contrary to the irreducibility of $f_3$. $\Box$ \end{pf}

These techniques are not sufficient to generalize to higher $k$.  In particular there are two problem cases with $k=4$, the case where $D_1 D_2 D_3 D_4$ is the only square and the case where $D_1 D_2 D_3 D_4$, $D_1 D_2$, and $D_3 D_4$ are the only squares.

\begin{rem}

Suppose $F(n)$ is the product of distinct irreducible quadratic polynomials $f_i$.  Roughly, we expect that the large primes factors of $\prod_{n \le x} f_i(n)$ should be rather sparse and should not overlap much with the large prime factors of $\prod_{n \le x} f_j(n)$.

By interpreting the DFIT result - incorrectly - as a statement of probability, one can refine this heuristic argument to estimate that the squarefree part of $N_x$ should tend towards $N_x^{1/2 + o(1)}$ as $x$ tends to infinity.  We cannot yet prove such a statement, and so leave it here as a conjecture.

\end{rem}

\section{Quadratic and Linear Terms}

Now suppose we wish to extend our results still farther, to consider products of the form $$ N_x = \prod_{n \le x} \left( \prod_{i=1}^I f_i(n) \right) \left( \prod_{k=1}^K g_k(n) \right)$$ where the $f_i$ are quadratic and the $g_k$ are linear.  Under what conditions for the $f_i,g_k$ will $N_x$ again be only finitely often a square?

We will assume, as we did with the $f_i$, that $g_k(n) \neq 0$ for any $n \ge 1$.

If the $f_i$ satisfy the conditions of Corollary \ref{ThmStrong} and $g_i(n)=n+b_i$ for all $i$ then the conclusion still holds, as under these conditions only a finite number (independent of $x$) of primes larger than $x$ could divide any of the terms $\prod_{n \le x} g_k(n)$.

\begin{thm}

Suppose that we have a set of $I$ distinct quadratic polynomials $f_i(n)$ with $f_1$ monic such that $J_f > 0.$ Suppose we also have a set of $K$ distinct linear polynomials $g_k(n)=a_k n + b_k$, where each $a_k\ge 2$ is relatively prime to each $D_i$ and to every other $a_k \ge 2$

Then the number $$ N_x = \prod_{n \le x} \left( \prod_{i=1}^I f_i(n) \right) \left( \prod_{i=1}^K g_i(n) \right)$$ cannot be infinitely often a squarefull.

\end{thm}

\begin{pf}

A prime $p > (2-\delta)x$ divides the term $g_k(n)$ when the following congruence holds
$$a_k n + b_k \equiv 0 (\textrm{mod} p)$$
$$n \equiv \frac{-b_k}{a_k} (\textrm{mod} p)$$

We can solve this explicitly since $n \le x < p$ implies that $n$ will equal $(kp-b_i)/a_i$ where $k$ is the smallest positive integer for which $jp-b_i$ is divisible by $a_i$.  However, we need $n \le x$, while $p > (2-\delta)x$ so this means that we would need to have $$\frac{j(2-\delta)x-b_i}{a_i} < \frac{jp-b_i}{a_i} \le x$$  $$ j-\frac{b_i}{(2-\delta)x} < \frac{a_i}{(2-\delta)}$$

Since $j$ is discrete, we can pick $x$ large enough so that the term $\frac{b_i}{2x}< 1/5$ and pick $\delta$ small enough so that $|a_i/(2-\delta) - a_i/2| < 1/5$ as well.  Then we get $$ j \le \frac{a_i}{2}+ \frac{2}{5}$$ or, in other words, that only half of the congruence classes modulo $a_i$ can contain primes larger than $(2-\delta)x$ which divide $\prod_{n \le x} g_i(n)$.

By Lemma \ref{congruence}, there must be a congruence class $k' \pmod{\prod D_i}$ for which primes congruent to $k' \pmod{\prod D_i}$ can, and eventually will, divide $\prod_{n\le x} f_1(n)$ but will never divide any other $\prod_{n \le x} f_i(n)$.  Provided $a_i$ is relatively prime to $\prod D_i$ the associated congruence classes modulo $a_i$ coming from $0 \le j \le \frac{a_i}{2}$ cannot cover the congruence class $k' \pmod{\prod D_i}$ completely.  So there exists some congruence class modulo $a_i D$ such that all sufficiently large primes in that congruence class eventually must divide $\prod_{n\le x} f_1(n)$ as $x$ grows, but which will never divide $\prod_{n \le x} g_k(n)$.

Since $a_j \neq a_i$ is relatively prime to $a_i \prod D_i$ we can repeat this process, and continue repeating through all of the $a_k$'s until we have found a congruence class $k'' \pmod{\prod a_k \prod D_i}$ such that all sufficiently large primes in that congruence class will eventually divide $\prod_{n\le x} f_1(n)$ but cannot divide $N_x / \prod_{n \le x} f_1(n)$.

Then, as in the proof of Theorem \ref{NotSquarefull}, for sufficiently large $x$, some of these primes can only divide $N_x$ precisely one time, and thus $N_x$ cannot be squarefull. $\Box$ \end{pf}

\begin{rem} We can, without difficulty, allow two linear terms with the same leading coefficient, say $a n+b ,a n+b'$ provided $a$ is prime (and as before relatively prime to all other $a_k$'s)  and $b \neq -b' \pmod{a}$.  This last condition will ensure that there is still some congruence class modulo $a$, such that primes from that congruence class can never divide $\prod_{n \le x} (an+b)(an+b')$. \end{rem}

Using slightly different techniques, we can prove the following theorem.

\begin{thm}

Let $f_i(n)$, $i \in \{1, 2, \ldots , I\}$, be distinct quadratic polynomials, and $g_k(n)=a_k n + b_k$, $ k \in \{1, 2, \ldots , K\}$, be distinct linear polynomials with non-zero, relatively prime coefficients, such that
\begin{enumerate}
\item $$J'_f := 1+ \sum_{\varnothing \neq J \subset \{1,2,3,\ldots,I\}, \prod_{j \in J} D_j \textrm{square}}  (-1)^{|J|} \neq 0$$
\item $a_1$ is positive and $a_1 \ge |a_k|$ for all $1 < k \le K$.
\item $a_1$ is relatively prime to $\prod_{i \le I} D_i$
\item For all $k$ such that $a_1 = |a_k|$, we have that $b_1 \neq b_k \pmod{a_1}$.
\end{enumerate}
Then $$ N_x = \prod_{n \le x} \left( \prod_{i=1}^I f_i(n) \right) \left( \prod_{k=1}^K g_k(n) \right)$$ can only be a perfect square finitely many times.

\end{thm}

\begin{pf}

Let us write $g_1(n)=an+b$.  Clearly all primes congruent to $b$ modulo $a$ less than $ax+b$ but larger than $a+b$ divide $\prod_{n \le x} g_1(n)$. Moreover, each prime of this congruence class that exists between $(a-\frac{1}{2})x$ and $ax+b$ divides $\prod_{n \le x} g_1(n)$ exactly once for sufficiently large $x$.  To see this, suppose $p$ is a prime congruent to $b \pmod{a}$ in the range $((a-\frac{1}{2})x,ax+b)$, and let $n'=(p-b)/a$.  Then clearly the first time $g_1(n)$ is divisible by $p$ is when $n = n'$.  The next time it happens is when $n = n' + p > p  \ge (a-\frac{1}{2})x > x$ which means $p$ divides $\prod_{n \le x} g_1(n)$ exactly one time.

In fact, these primes can only divide $\prod_{n \le x} \left( \prod_{k=1}^K g_k(n) \right)$ once for large enough $x$.  Every $g_k$ with $|a_k| < a_1$ can only contribute primes smaller than $(a-\frac{1}{2})x$.  By assumption, every $g_k$ with the same leading coefficient as $g_1$ is of the form $an+b'$ where $b' \neq b \pmod{a}$.  Thus the first time a prime $p>(a-\frac{1}{2})x$ congruent to $b \pmod{a}$ divides $g_k(n)$ is at the earliest when $n = (2p-b')/a > ((2a-1)x-b')/a >x $ once $x$ is large enough.

So in order for $N_x$ to be a square, each of the primes congruent to $b \pmod{a}$ in the range $((a-\frac{1}{2})x,ax+b)$ must divide $\prod_{n \le x} \left( \prod_{i=1}^I f_i(n) \right)$.  However, by a similar argument to Lemma \label{congruence}, the proportion of the primes that can never divide any of these terms is $$\left| \sum_{D < p \le z} \prod_{1 \le i \le I} (-1) \frac{-1+\left( \frac{D_i}{p} \right)}{2} \right|,$$ which will be asymptotic to a non-zero proportion of $\pi(z)$ whenever $J'_f \neq 0$.

These correspond to a proportion of residue classes modulo $\prod_{i\le I} D_i$.  Since $a$ is relatively prime to $\prod_{i\le I} D_i$, there must exist residue classes modulo $a \prod_{i\le I} D_i$ such that they reduce to $b \pmod{a}$ and yet primes in these residue classes can never divide $\prod_{n \le x} \left( \prod_{i=1}^I f_i(n) \right)$.  Now, if we pick $x$ to be large enough, then there must exist a prime congruent to $b$ modulo $a$ in the region $((a-1/2)x,ax+b)$, which cannot divide $\prod_{n \le x} \left( \prod_{i=1}^I f_i(n) \right)$ yet must divide $\prod_{n \le x} \left( \prod_{k=1}^K g_k(n) \right)$ precisely once.

Thus $N_x$ cannot be infinitely often a square. $\Box$ \end{pf}

\section{More general F(n)}

In the case of still more general $F(n)$ we cannot yet obtain any theorems which say that $N_x$ will only be finitely often a square or finitely often squarefull, yet we can obtain a small density result.

Here, given a function $F(n) \in \mathbb{Z}[n]$, let $d_F$ be the positive integer such that there is some element of the Galois group of $F$ which fixes precisely $d_F$ roots of $F(n)$ and any element which fixes less than $d_F$ roots of $F(n)$ will fix none of the roots.  $d_F$ exists since the Galois group contains the trivial element which will fix all the roots of $F$, which also implies that $d_F \le \deg{F}$.

We denote the size of the Galois group of $F$ by $g_F$.

\begin{thm}\label{pth}

Suppose $F(n) \in \mathbb{Z}[n]$ is not of the form $s (G(n))^p$ for some rational number $s$ and some polynomial $G(n) \in \mathbb{Z}[n]$.  Then $$\# \{x \le X | N_x \text{is a perfect $p^{th}$ power}\} = O\left( X^{\log{(d_F+1)}/\log{\lceil p/d_F \rceil}} \right)$$ and, more generally, $$\# \{x \le X | N_x \text{is a perfect $p^{th}$ power}\} = O\left( X^{24/25}\right)$$

\end{thm}

We note that this generalizes the results of Cilleruelo, et al., in \cite{Cil2}.

To begin, we need the following lemma.

\begin{lem}\label{primeseq}

There exists a sequence of primes $q_1, q_2, q_3, \ldots$, such that $$\lceil p/d \rceil q_i \left(1-g_F \frac{\log{q_i}}{q_i}\right) \le q_{i+1} \le \lceil p/d_F \rceil q_i$$ and $F(n)$ has $d_F$ roots modulo $q_i$.

\end{lem}

\begin{pf} If $f(n)$ is some irreducible polynomial, then the way that $f$ factors when taken modulo some prime $p$ is determined completely by the way the Frobenius automorphism acts on the roots of $f$.  Taken modulo $p$, the Frobenius automorphism maps the set of roots of $f$ bijectively onto the roots of $f$.  If it maps an element onto itself, this corresponds to a linear factor of $f$ modulo $p$.  If it maps one element onto a second element, and the second element back onto the first (i.e. a 2-cycle), then this corresponds to a quadratic factor of $f$ modulo $p$, and so on.

Thus if the cycle structure of the Frobenius automorphism acting on the roots of $f$ is $(m_1,m_2,\ldots,m_r)$, then $$f(n) \equiv \prod_{i=1}^r g_i(n) \pmod{p}$$ where $\deg{g_i} = m_i$ and each $g_i$ is irreducible modulo $p$.

A similar result holds even if our function is reducible.  In particular, let $F(n) = \prod f_i(n)^{e_i}$ for distinct irreducibles $f_i$.  We can still consider the Galois group of $F$ as the compositum of all the Galois groups of the $f_i$'s; this is also the splitting field for $F$.  The Frobenius automorphism for a given prime $p$ is again an element of the Galois group of $F$ and it will map roots of $F$ bijectively onto roots of $F$, and will actually map roots of $f_i$ bijectively onto roots of $f_i$.  Thus if the cycle structure of the Frobenius automorphism acting on the roots of $F$ is $(m_1,m_2,\ldots,m_r)$, then $$F(n) \equiv \prod_{i=1}^r g_i(n) \pmod{p}$$ where $\deg{g_i} = m_i$ and each $g_i$ is irreducible modulo $p$.

In particular, this tells us that $F(n)$ has $d$ roots modulo $p$ precisely when the Frobenius automorphism fixes exactly $d$ of the roots of $F(n)$.

The Chebotarev Density theorem (see \cite{IK}, page 143) says that there exists a natural density of primes $p$ for which the cycle structure of the Frobenius Automorphism of $p$ acting on the roots of $F$ is $(m_1,m_2,\ldots,m_r)$.  In particular this density is the number of elements of the Galois group which produce this cycle structure when they act on the roots of $F$ divided by the total number of elements in the Galois group.

Since we know, by definition, that there exists some element of the Galois group that fixes precisely $d_F$ roots of $F$, there must be a positive density $c$ of primes $p$ for which $F$ has precisely $d_F$ roots modulo $p$.

Thus, if we let $\epsilon(X)=g_F \log{X}/X$ and let $\pi_{d_F} (X)$ denote the number of primes less than $X$ for which $F$ has $d_F$ roots, then \begin{align*}
&\pi_{d_F} (X) -\pi_{d_F} (X(1-\epsilon(X))) \\
&\sim c \frac{X}{\log{X}} - c \frac{X(1-\epsilon(X))}{\log{X}} \\
&=c\epsilon(X) \frac{X}{\log{X}}\\
&=cg_F \ge 1
\end{align*} since $c \ge 1/g_F$.

Thus we can find a prime $q_{i+1}$ which is between $\lceil p/d_F \rceil q_i(1 - g_F \log{q_i}/q_i)$ and $\lceil p/d \rceil q_i$ and for which $F$ has $d_F$ roots modulo $q_{i+1}$, provided we start this sequence with a sufficiently large prime $q_1$.  $\Box$ \end{pf}

Here, if $F(n)=s f_1(n)^{e_1} \cdots f_k(n)^{e_k}$ for some $s \in \mathbb{Q}$ and for distinct irreducible polynomials $f_i$, we let $\text{sdisc}(F)$ denote the discriminant of $\prod_{i=1}^k f_i(n)$.

Recall that if $F(n)$ has $k$ roots modulo $p$, then it also has $k$ roots modulo $p^i$ provided $p$ does not divide $\text{sdisc}(F)$.  This is true because if $p$ does not divide the discriminant of $f_i$ then the roots of $f_i$ modulo $p$ are distinct, and we can then apply Hensel's lemma to see that these roots extend to distinct roots modulo $p^i$.

Now consider any of the primes $q_i$.  Let $a_i(x)$ represent the number of times $q_i$ divides $N_x$.

By our construction of the $q_i$, we know that $F(n)$ has $d_F$ roots modulo $q_i$.  Thus, $a_i(x+q_i)-a_i(x) \ge d_F$.

At the same time we know that $F(n)$ has $d_F$ roots modulo $q_i^2$, so $a_i(x+j+1)-a_i(x+j)>1$ for at most $d_F$ values of $j$, with $0\le j \le q_i^2-1$.

Let us further assume that if $p|a_i(x)$, then $x$ belongs to an interval of length $q_i$ on which $a_i$ is constant, and suppose these intervals are distinct; this will overestimate how often $p|a_i(x)$ but still give us our big-Oh bounds.  Now, we will estimate how close two successive intervals can be on average.  Let $I_1$, $I_2$ be the two intervals in question, with $I_1 = [x_1,x_1+q_i-1]$ and $I_2 = [x_2,x_2+q_i-1]$.  If for all $x_1 \le x < x_2$, we have that $a_i(x+1)-a_i(x)\le 1$, so then for all $x_1 \le x \le x_2 -q_i +1$ we have that $a_i(x+q_i)-a_i(x) = d_F$.  Thus $a_i(x_2)-a_i(x_1)\le d_F \lceil (x_2-x_1)/q_i \rceil$ and at the same time $a_i(x_2) - a_i(x_1)= p$.  Thus, $x_2-x_1 \ge \lfloor p/d_F \rfloor q_i$.

However, we also know that over an interval of $x$'s of length $q_i^2$, $a_i$ will jump by more than one exactly $d_F$ times.  Thus it is possible that we could have two sub-intervals $I_1 = [x_1,x_1+q_i-1]$ and $I_2 = [x_2,x_2+q_i-1]$ of the type discussed in the previous paragraph with $x_2 = x_1 + q_i$, but this could only occur at most $d_F$ times over the full interval.  Each other pair of successive intervals must be separated as in the previous paragraph.

Thus we see that if we have an interval of length $q_{i+1}$, which is slightly smaller than $\lceil p/d_F \rceil q_i$, then it can contain at most $d_F+1$ sub-intervals of length $q_i$ of values of $x$ for which $p|a_i(x)$; consequently, if $X>q_{i+1}$ at most $$2X \frac{(d_F+1)q_i}{q_{i+1}}$$ of the numbers $x$ up to $X$, will have $N_x$ be a perfect $p^{th}$ power.  (Here the $2$ is a fudge factor since $X$ will likely not be a multiple of $q_{i+1}$.)

If we look at an interval of length $q_{i+2}$ then it can contain at most $d_F+1$ intervals of length $q_{i+1}$ of values of $x$ for which $p|a_{i+1}(x)$, which themselves can contain at most $d_F+1$ intervals of length $q_i$ of values of $x$ for which $p|a_i(x)$; consequently, at most $$2X \frac{(d_F+1)q_i}{q_{i+1}} \frac{(d_F+1)q_{i+1}}{q_{i+2}}$$ of the numbers $x$ up to $X$, if $X>q_{i+2}$, will have $N_x$ be a perfect $p^{th}$ power.  And so on.

Now suppose $q_i \le X < q_{i+1}$ then we have that there are at most $$2 X \left(\frac{d_F+1}{\lceil p/d \rceil}\right)^{i-1} (1-\epsilon(q_1))^{-1}(1-\epsilon(q_2))^{-1}\cdots (1-\epsilon(q_{i-1}))^{-1}$$ $x$ less than $X$ for which $N_x$ is a perfect $p^{th}$ power.

Note that $i-1 > \log_{\lceil p/d_F \rceil}{(X/q_1)}$, so \begin{align*}
&\left(\frac{d_f+1}{\lceil p/d_F \rceil}\right)^{i-1} \\
&< \left(\frac{d_F+1}{\lceil p/d_F \rceil}\right)^{\log_{\lceil p/d_F \rceil}{(X/q_1)}} \\
&= \exp{\left( \frac{\log{X}-\log{q_1}}{\log{\lceil p/d_F \rceil}} \log{\frac{d_F+1}{\lceil p/d_F \rceil}} \right)} \\
&= X^{(\log{(d_F+1)}/\log{\lceil p/d_F \rceil} -1)} \left(\frac{\lceil p/d_F \rceil}{d_F+1} \right)^{\log{q_1}/\log{\lceil p/d_F \rceil}}
\end{align*}

Furthermore, note that \begin{align*}
&\epsilon(q_i)= \frac{g_F\log{q_i}}{q_i} \\
&= O\left(\frac{g_F\log{(q_1\lceil p/d_F \rceil^{i-1})}}{(q_1\lfloor p/d_F \rfloor^{i-1})}\right) \\
&= O\left(\frac{i}{\lfloor p/d_F \rfloor^{(i-1)}}\right)
\end{align*}

Thus \begin{align*}
&(1-\epsilon(q_1))^{-1}(1-\epsilon(q_2))^{-1}\cdots (1-\epsilon(q_{i-1}))^{-1} \\
&\le \exp{\left(\sum_{n=1}^{i-1}  \epsilon(q_i) \right)} = \exp{\left(O\left(\sum_{n=1}^{\infty}  \frac{n}{\lfloor p/d_F \rfloor^{n-1}} \right)\right)} \\
&=\exp{\left(O\left(\left(\sum_{n=1}^{\infty}  \frac{1}{\lfloor p/d_F \rfloor^{n-1}}\right)^2 \right)\right)} = \exp{\left(O\left(\left(\frac{1}{1-\frac{1}{\lfloor p/d_F \rfloor}}\right)^2 \right)\right)}
\end{align*} which is clearly bounded.

Together these estimates prove the first part of Theorem \ref{pth}; however this result is only interesting when $p^2 > d_F$, for smaller $p$ we will use a variation of the Tur\'{a}n sieve (following the method of \cite{Coj}).

For the Tur\'{a}n sieve, let $\mathcal{A}$ be an arbitrary finite set, $\mathcal{P}$ be some set of primes, and to each prime $p \in \mathcal{P}$ associate a set $\mathcal{A}_p \subset \mathcal{A}$, and let $\mathcal{A}_{p,q} \: = \mathcal{A}_p \cap \mathcal{A}_q$.  Now suppose $$\# \mathcal{A}_p = \delta_p X+ R_p$$ and $$\# \mathcal{A}_{p,q}= \delta_p \delta_q X+ R_{p,q}$$ where $X= \# \mathcal{A}$, then we have the following result.

\begin{thm}

With all notation as in the previous paragraph, let $$U(z)\: = \sum_{\substack{p \in \mathcal{P} \\ p \le z}} \delta_p$$ then $$\# \left( \mathcal{A} \setminus \bigcup_{p \in \mathcal{P}} \mathcal{A}_p \right) \le \frac{X}{U(z)} + \frac{2}{U(z)} \sum_{\substack{p \in \mathcal{P} \\ p \le z}} \left| R_p \right| + \frac{1}{U(z)^2} \sum_{\substack{p,q \in \mathcal{P} \\ p,q \le z}} \left| R_{p,q} \right|$$

\end{thm}

We also need the following result.  Here we use the shorthand $F_k(n):=F(n)F(n+1)\cdots F(n+k)$.

\begin{lem}

Suppose $F(n) \in \mathbb{Z}[n]$ is not of the form $s G(n)^p$ for some $s \in \mathbb{Q}$ and $G(n) \in \mathbb{Z}[n]$.  Then for any prime $q$ which does not divide $\text{sdisc}(F)$ and is larger than $\deg{(F)} k$, we have that $F_k(n)$ taken modulo $q$ is not equivalent to $s G(n)^p$ for any $G(n) \in \mathbb{Z}_q [n]$, $s \in \mathbb{Z}_q$.

\end{lem}

\begin{pf}

If $q \nmid \text{sdisc}(F)$, then the factors of $f_i$ modulo $q$ must be distinct from the factors of $f_j$ modulo $q$ if $i \neq j$, and the factors of $f_i$ modulo $q$ are themselves distinct from each other.  Thus if not every $f_i$ divides $F(n)$ with a $p$-multiple multiplicity, then not every irreducible modulo $q$ divides $F(n)$ with a $p$-multiple multiplicity.

Moreover, $g(n)$ is irreducible (over $\mathbb{Z}$ or $\mathbb{Z}_q$) if and only if $g(n+1)$ is also irreducible, and $g(n)^m|F(n)$ if and only if $g(n+1)^m | F(n+1)$.  If $g(n)=n^l+a_{l-1} n^{l-1} + \cdots + a_0$ then $g(n+i) = n^l + (il+ a_{l-1}) n^{l-1} + \cdots$ so these will be distinct modulo $q$ if $il \neq 0 \pmod{q}$.

Now, consider $F_k(n)$ and suppose that all irreducibles divide $F_k(n)$ with a $p$-multiple multiplicity when we reduce $F_k(n)$ modulo $q$.  By the work above we know that there exists some irreducible over $\mathbb{Z}_q [n]$, let us call it $g_1(n)$, that does not divide $F(n)$ with a $p$-multiple multiplicity, but it does divide $F_k(n)$ with a $p$-multiple multiplicity, therefore there must be some other irreducible $g_2(n)$ such that $g_2(n)|F(n)$ and $g_2(n+i)\equiv g_1(n) \pmod{q}$ with $1 \le i \le k$; however, we can assume that $g_2(n)$ does not divide $F(n)$ a multiple of $p$ times, so we can find a $g_3,g_4,\ldots$ in this fashion.  However, $F(n)$ has finite degree, so this sequence of $g_i$'s must eventually repeat itself.  Suppose, without loss of generality, that $g_1(n) = g_j(n)$ with $j$ minimal, then we have that $g_1(n + i) \equiv g_j(n) \pmod{q}$ for $j-1 \le i \le k(j-1)$.  By the previous paragraph, that implies $i \equiv 0 \pmod{q}$ but $i > 0 $ and $i \le k(j-1) < \deg{(F)} k$, since $F$ can have at most $\deg{(F)}$ irreducible factors.  So since $q > \deg{(F)} k$, $F_k(n)$ cannot be of the form $s G(n)^p$ for any $G(n) \in \mathbb{Z}_q [n]$, $s \in \mathbb{Z}_q$, as desired. $\Box$ \end{pf}

In our case, let $$\mathcal{A}\: = \{n \le X\},$$ and for each prime $q \nmid \text{sdisc}(F_k)$, let $$\mathcal{A}_q \: = \{n \le X | F_k(n) \text{is not a perfect $p^{th}$ power modulo $q$}\}.$$  Note that $\text{sdisc}(F_k) = \text{sdisc}(F)$, since $\text{disc}(f_i(n))=\text{disc}(f_i(n+1))$.

By \cite{Li}, page 94, we have that $$\left| \sum_{a \mod{q}} \chi_p \left( F_k(a) \right) \right| \le (\deg{F_k}-1)\sqrt{q}$$ if $\chi_p$ is a non-trivial multiplicative character of order $p$ and $F_k$ has some root modulo $q$ whose multiplicity is not a $p$-multiple.  By the previous lemma, this latter requirement is satisfied.

Let $S_k$ denote the number of $n$ modulo $q$ for which $F_k(n)$ is $p^{th}$ power modulo $q$, then supposing there exist non-trivial characters, we have \begin{align*} &|p S_k - q| = \left|\sum_{\chi_p \neq 1} \sum_{a\pmod{q}} \chi_p \left( \frac{F_k(a)}{q} \right)\right| \\
&\le (p-1)(\deg{F_k}-1)\sqrt{q}
\end{align*}

Thus $S_k = q/p+O((\deg{F_k})\sqrt{q})$.

Thus \begin{align*}\# \mathcal{A}_q &= \left(\frac{X}{q}+O(1)\right) \left( \frac{q(p-1)}{p} + O((\deg{F_k})\sqrt{q})\right) \\
&= \frac{X(p-1)}{p}+O\left(\deg{F_k}\frac{X}{\sqrt{q}} + q\deg{F_k}\right)\end{align*} and similarly, given distinct primes $q_1,q_2$ we have \begin{align*} \# \mathcal{A}_{q_1,q_2} &= \left(\frac{X}{q_1 q_2}+O(1)\right) \left( \frac{q_1 q_2 (p-1)^2}{p^2} +  O((\deg{F_k})^2(\sqrt{q_1}q_2 + \sqrt{q_2}q_1))\right)\\
&= \frac{X(p-1)^2}{p^2} +  O\left((\deg{F_k}^2)\left(\frac{X}{\sqrt{q_1}} + \frac{X}{\sqrt{q_2}} + q_1 q_2 \right)\right) \end{align*}

For our set of primes $\mathcal{P}$ we want the set of all primes $q$ between $z$ and $2z$, such that $q$ does not divide $\text{sdisc}(F_k)$ and $q \equiv 1 \pmod{p}$ (so that there will exist non-trivial characters).  We will determine $z$ later.

Then by the Tur\'{a}n sieve, the number of $n \le X$ for which $F_k(n)$ is a perfect $p^{th}$ power is $$\ll X\frac{\log{z}}{z} + (\deg{F_k})\frac{X}{\sqrt{z}} + (\deg{F_k})z + (\deg{F_k})^2\frac{X}{\sqrt{z}} + (\deg{F_k})^2 z^2$$ and the implied constant is independent of our choice for $k$.

We now use the following lemma to see how frequently $N_x, N_{x+k}$ can be both a $p^{th}$ power, with $k$ small.

\begin{lem}

Let $S(X)$ be some subset of the natural numbers $\{1,2,\ldots,X\}$, and suppose $|S(X)|> X/K(X)$ for some function $K(X)< X$.

Let $S(X)_k$ denote those $s \in S(X)$ such that $s+k \in S(X)$ and $s+i \in X \setminus S(X)$ for $1 \le i < k$.

Then there exists some integer $k \le K(X)$ such that $\left| S(X)_k \right| \ge 2X/K(X)^3.$

\end{lem}

\begin{pf}

Suppose to the contrary that for all $k \le K(X)$ there are less than $2X/K(X)^3$ elements in $S(X)_k$.  Let us consider the most number of elements that could be in $S(X)$ under these conditions.  In particular we want to have as small a gap between successive elements as possible.  So let us assume that for all $k \le K(X)$ there are at most $2X/K(X)^3-1 $ distinct $s \in S(X)$ for which $s+k \in S(X)$ and $s+i \in X \setminus S(X)$ for $1 \le i < k$.  The number of integers in the union $$\bigcup_{k \le K(X)} \bigcup_{s\in S(X)_k} \{s,s+1,\ldots,s+k-1\}$$ is then at most $$\frac{K(X)(K(X)+1)}{2} \lceil \frac{2X}{K(X)^3}-1 \rceil)= \frac{X(K(X)+1)}{K(X)^2}-\frac{K(X)(K(X)+1)}{2} $$

Then let us also suppose, in order to maximize the number of elements in $S(X)$, that for each remaining $s \in S(X)$, the first element in $S(X)$ after $s$ is $s+K(X)+1$.

Thus the total number of elements in $S(X)$ is, at most, \begin{align*} &K(X) \left( \frac{2X}{K(X)^3}-1 \right) + (X-\frac{X(K(X)+1)}{K(X)^2}+\frac{K(X)(K(X)+1)}{2})\frac{1}{K(X)+1}+1\\
&= \frac{2X}{K(X)^2} - K(X) + \frac{X}{K(X)+1} -\frac{X}{K(X)^2} + \frac{K(X)}{2}+1\\
&=\frac{X}{K(X)^2} - \frac{K(X)}{2} + \frac{X}{K(X)+1}+1
\end{align*} which is smaller than $X/K(X)$, since \begin{align*} \frac{X}{K(X)}-\frac{X}{K(X)+1} &=\frac{X}{K(X)}\left( 1-\frac{1}{1+\frac{1}{K(X)}}\right) \\
&=\frac{X}{K(X)}\left( \frac{1}{K(X)}-\frac{1}{K(X)^2}+\cdots\right) \\
&<\frac{X}{K(X)^2} \end{align*}

$\Box$

\end{pf}

Now we consider $F(n)$ again.  Suppose that $N_x$ is a perfect $p^{th}$ power for at least $X/K(x)$ of the $x \le X$.  Then the lemma above implies that there must be some $k < K(X)$ for which there are at least $2X/K(X)^3$  of the $x \le X$ such that $N_x,N_{x+k}$ are both perfect $p^{th}$ powers and there are no such powers between them.  Since $N_x,N_{x+k}$ are both perfect $p^{th}$ powers, so must $F_k(x)=F(x+1)F(x+2)\ldots F(x+k)$ be a perfect $p^{th}$ power.

According to the above work $F_k(n)$ is a perfect $p^{th}$ power $$\ll X\frac{\log{z}}{z} + (\deg{F_k})\frac{X}{\sqrt{z}} + (\deg{F_k})z + (\deg{F_k})^2\frac{X}{\sqrt{z}} + (\deg{F_k})^2 z^2$$ times which is \begin{align*}&\ll_F X\frac{\log{z}}{z} + k\frac{X}{\sqrt{z}} + kz + k^2\frac{X}{\sqrt{z}} +k^2 z^2 \\
&\ll_F K(X)^2\frac{X}{\sqrt{z}} +K(X)^2 z^2\end{align*}
but by assumption $F_k(n)$ is a perfect $p^{th}$ power at least $2X/K(X)^3$ times.  Putting these together we see that $K(X)$ must satisfy $$X \ll  K(X)^5 \frac{X}{\sqrt{z}} +K(X)^5 z^2$$ for any choice of $z$.

Setting $z=X^{2/5}$ we see that $K(X)$ cannot have smaller magnitude than $$X^{1/25}.$$  Thus we have proved the second part of Theorem \ref{pth}.

\section{Acknowledgements}

The authors acknowledge support from National Science Foundation grant DMS 08-38434 ``EMSW21-MCTP: Research Experience for Graduate Students''.

The authors also wish to thank Paul Pollack for his numerous suggestions and assistance.

\section{References}
\bibliographystyle{elsarticle-num}
\bibliography{GoodReferences2}

\end{document}